\documentclass[11pt,amsfonts]{article}
\usepackage{graphicx}
\usepackage{latexsym}
\usepackage{amssymb}
\usepackage{amsmath}
\usepackage{layout}
\newtheorem{prop}{Proposition}
\newtheorem{lemma}{Lemma}
\newtheorem{definition}{Definition}
\newtheorem{corollary}{Corollary}
\newtheorem{theorem}{Theorem}
\newtheorem{remark}{Remark}

\def\real{{\mathord{{\rm I\kern-2.8pt R}}}}        
\def\inte{{\mathord{{\rm I\kern-2.8pt N}}}}

\def\sZZ{{\rm Z\kern-2.8ptem{}Z}}

\def\z{{\mathchoice
  {\sZZ}
  {\sZZ}
  {\rm Z\kern-0.30em{}Z}
  {\rm Z\kern-0.25em{}Z} }}
\def\sQQ{{\kern 0.27em \vrule height1.45ex width0.03em depth0em
          \kern-0.30em \rm Q}}
\def\qu{{\mathchoice
    {\sQQ}
    {\sQQ}
  {\kern 0.225em \vrule height1.05ex width0.025em depth0em \kern-0.25em \rm Q}
  {\kern 0.180em \vrule height0.78ex width0.020em depth0em \kern-0.20em \rm Q}
        }}
\def\sCC{{\kern 0.27em \vrule height1.45ex width0.03em depth0em
          \kern-0.30em \rm C}}
\def\complex{{\mathchoice
    {\sCC}
    {\sCC}
  {\kern 0.225em \vrule height1.05ex width0.025em depth0em \kern-0.25em \rm C}
  {\kern 0.180em \vrule height0.78ex width0.020em depth0em \kern-0.20em \rm C}
        }}

\newcommand{\ba}{\begin{array}}
\newcommand{\ea}{\end{array}}
\newcommand{\be}{\begin{equation}}
\newcommand{\ee}{\end{equation}}
\newcommand{\bea}{\begin{eqnarray}}
\newcommand{\eea}{\end{eqnarray}}
\newcommand{\beaa}{\begin{eqnarray*}}
\newcommand{\eeaa}{\end{eqnarray*}}

\def\z{\zeta}

\font\tenmath=msbm10 \font\sevenmath=msbm7 \font\fivemath=msbm5
\newfam\mathfam \textfont\mathfam=\tenmath
\scriptfont\mathfam=\sevenmath \scriptscriptfont\mathfam=\fivemath

\def \={{\buildrel {\rm (law)} \over =}}

\def\qed{ \hfill \vrule width.25cm height.25cm depth0cm\smallskip}

\newcommand{\basa}{\begin{assumption}}
\newcommand{\easa}{\end{assumption}}

\newcommand{\bas}{\begin{assum}}
\newcommand{\eas}{\end{assum}}

\def\liminf{\mathop{\underline{\rm lim}}}

\newcommand{\ignore}[1]{}
\begin{document}

\date{ }
\title{Asymptotic Cram\'er's theorem and analysis on Wiener space}
\author{ Ciprian A. TUDOR \\
 Laboratoire Paul Painlev\'e, Universit\'e de Lille 1,\\
  F-59655 Villeneuve d'Ascq, France,\\
   email: tudor@math.univ-lille1.fr.}
\maketitle

\begin{abstract}
We prove an asymptotic Cram\'er's theorem, that is, if the sequence $(X_{n}+ Y_{n})_{n\geq 1}$ converges in law to the standard normal distribution and for every $n\geq 1$ the random variables $X_{n}$ and $Y_{n}$ are independent, then $(X_{n})_{n\geq 1}$ {\it and }  $(Y_{n}) _{n\geq 1}$ converge in law to a normal distribution. Then we compare this result with  recent criteria for the central convergence obtained in terms of Malliavin derivatives.
\end{abstract}

{\bf  2000 AMS Classification Numbers:} 60G15,  60H05, 60F05, 60H07.


\section{Introduction}
The sum of two independent random variables with Gaussian distribution is a  Gaussian random variable.  A famous result by Harald Cram\'er \cite{Cramer} says that the converse implication is also true. Namely, if the law of  $X+Y$ is Gaussian and $X$ and $Y$ are independent random variables, then $X$ and $Y$ are Gaussian.
We study in this paper  the following problem: given two sequences of centered square integrable random variables $(X_{n})_{n\geq 1}$ and $(Y_{n})_{n\geq 1} $ such that $\mathbf{E}X_{n}^{2} \to _{n\to \infty } c_{1} $ and $\mathbf{E}Y_{n}^{2} \to _{n\to \infty} c_{2}$ with $c_{1}, c_{2}>0$ and $c_{1}+ c_{2}=1$ and assuming that for every $n\geq 1$, $X_{n}$ and $Y_{n}$ are independent and $X_{n}+ Y_{n}\to _{n\to \infty }N(0,1) $ in law, can we get the convergence of $X_{n} $ to the normal law $N(0, c_{1})$ {\it and } the convergence of $Y_{n}$ to the normal law $N(0, c_{2})$? We will say in this case that the central limit of the sum is decoupled. A partial answer has been given in \cite{PT}: here the authors proved that the central limit for the sum  implies the central limit for each term when the random variables $X_{n}$ and $Y_{n}$ lives in a Wiener chaos of fixed order. In this work we will  prove this result for a very general class of random variables.

Then we will try to understand this asymptotic Cram\'er's theorem from the  perspective of some   recent ideas from \cite{NoPe1} and \cite{NoPe2} related to the Stein's method on Wiener space and some older ideas from \cite{NUZ}, \cite{UZ}, \cite{UZ2} where the independence of random variables is characterized in terms of the Malliavin derivatives. Let $\left( \Omega, {\cal{F}},P\right)$ be a probability space and let $(W_{t})_{t\in [0,1]}$ be  a Wiener process on this space.  Recall that a result in \cite{NoPe1} says that a sequence of Malliavin differentiable (with respect to $W$)  random variables $X_{n}$ (defined on $\Omega$)  converges to the normal law $N(0,1)$ if and only if
\begin{equation*}
\mathbf{E}\left( f'_{z} (X_{n}) (1-\langle DX_{n} , D(-L) ^{-1} X_{n} \rangle ) \right) \to _{n \to \infty }0.
\end{equation*}
where we denoted by $D$ the Malliavin derivative with respect to $W$, by $L$ the Ornstein-Uhlenbeck generator  and by $f_{z} $  the solution of the Stein's equation (for fixed $z\in \mathbb{R}$)
\begin{equation}\label{intro1}
1_{(-\infty , z] }(x)-P(Z\leq z) = f'(x)-xf(x), \hskip0.5cm x\in \mathbb{R}.
\end{equation}
(Throughout this paper we denote by $\langle \cdot , \cdot \rangle $ the scalar product in $L^{2}([0,1])$.) In particular, if $\mathbf{E}\left( 1-\langle DX_{n}, D(-L)^{-1} X_{n}\rangle \right) ^{2}\to _{n\to \infty }0$ then $X_{n}$ converges to $N(0,1)$ as $n\to \infty$ by using Schwarz's inequality and the fact that $f'_{z}$ is bounded (actually it suffices to have \\
$\mathbf{E}\left| 1-\langle DX_{n}, D(-L)^{-1} X_{n}\rangle \right| \to _{n\to \infty }0$).

Let us describe the basic idea  to treat the convergence of sums of independent random variables to the normal law. Let $X_{n}, Y_{n}$ be two sequences as above (that means Malliavin differentiable with $\mathbf{E}X_{n}^{2}\to_{n \to \infty} c_{1}>0$, $\mathbf{E}Y_{n}^{2} \to_{n \to \infty}c_{2}>0$ and $c_{1}+c_{2}=1$). The fact that $X_{n}+ Y_{n} \to _{n\to \infty } N(0,1)$ (in law) implies that
 \begin{equation}\label{int1}
 \mathbf{E}\left( f'_{z} (X_{n}+ Y_{n}) (1-\langle D(X_{n}+ Y_{n} ), D(-L) ^{-1} (X_{n}+ Y_{n}) \rangle ) \right) \to _{n \to \infty }0.
 \end{equation}
 Suppose now that $X_{n}$ and $Y_{n}$ are independent for every $n$. A result by \"Ustunel and Zakai (\cite{UZ}, Theorem 3) says that in this case
 \begin{equation}\label{intro2}
 \mathbf{E}(\langle DX_{n}, D(-L) ^{-1}Y_{n}\rangle | X_{n} )=0  \mbox{ and } \mathbf{E}(\langle DY_{n}, D(-L) ^{-1}X_{n}\rangle | Y_{n} )=0 \mbox{ a.s }.
\end{equation}
The relation (\ref{intro2}) induces the idea that the summands containing $\langle DX_{n} , D(-L) ^{-1}Y_{n} \rangle $ and \\
$\langle DY_{n} , D(-L) ^{-1}X_{n} \rangle $  could be eliminated from (\ref{int1}). We will see that it is  not immediate and  that actually a stronger condition than the independence of $X_{n}$ and $Y_{n}$ is necessary in order to do this. Therefore our first step is to introduce some classes of independent random variables $X,Y$ such that the "mixed" terms $\langle DX , D(-L) ^{-1}Y \rangle $ and $\langle DY , D(-L) ^{-1}X \rangle $ vanish. A first class that we consider here is the class of so-called {\it strongly independent} random variables for which every multiple integral in the chaos decomposition of $X$ is independent of every multiple integral in the chaos decomposition of $Y$. We will see that  if $X$ and $Y$ are strongly independent, then
$\langle DX , D(-L) ^{-1}Y \rangle =\langle DY , D(-L) ^{-1}X \rangle =0 $  almost surely.
Another class we consider is the class of random variables $X,Y$ that are differentiable in the Malliavin sense and such that $X$ is independent of the couple $(Y,\langle DY, D(-L) ^{-1}Y\rangle )$ and $Y$ is independent of the couple $(X, \langle DX, D(-L) ^{-1} X\rangle )$. We will say in this case that the couple $(X,Y)$ belongs to the class ${\cal{A}}$. A couple of strongly independent random variables belongs to ${\cal{A}}$ and in this sense this class is an intermediary class between the independent and strongly independent random variables. For couples in ${\cal{A}}$ we will show that $\mathbf{E}(\langle DX, D(-L)^{-1}Y\rangle | X+Y)=\mathbf{E}(\langle DY, D(-L)^{-1}X\rangle | X+Y) =0$ almost surely and it is then again possible to cancel the mixed terms in (\ref{int1}).
 We will prove, by an elementary argument coming from the original Cram\'er's theorem and  without using Malliavin calculus, that for independent random variables the asymptotic Cram\'er's theorem  holds. But  for random variables in these classes (in the class ${\cal{A}}$ or strongly independent) we can give  further results by using the tools of the Malliavin calculus. Concretely, we will treat the following problem: suppose that  the sum $X_{n}+Y_{n}$ converges to the normal law in a strong sense, that is, the upper bound $\mathbf{E}\left( 1-\langle D(X_{n}+ Y_{n}), D(-L) ^{-1} (X_{n}+ Y_{n}) \rangle \right) ^{2} $ converges to zero as $n\to \infty$. We can interpret this by saying that the sum $X_{n}+ Y_{n}$ is  "close" to $N(0,1)$, not in the sense of the rate of convergence but in the sense that $X_{n}+Y_{n}$ belongs to a subset of the set of the sequences of random variables converging to $N(0,1)$.  Then can we obtain the strong convergence of $X_{n}$ and $Y_{n}$ to the normal law, that is $\mathbf{E}\left( c_{1}-G_{X_{n}}\right) ^{2} \to _{n\to \infty}0$ and $\mathbf{E}\left( c_{2}-G_{Y_{n}}\right) ^{2}\to _{n\to \infty} 0$, where $G_{X_{n}}$ is given by (\ref{gx})? We prove that this property is true for strongly independent random variables while for couples in the class ${\cal{A}}$ a supplementary assumption is necessary in order to ensure the strong convergence of $X_{n}$ and $Y_{n}$ from the  convergence of $X_{n}+Y_{n}$.

 The organization of the paper is as follows. Section 2 contains preliminaries on the stochastic calculus of variations. In Section 3 we prove the asymptotic Cram\'er's theorem by using an elementary argument while  Section 4 contains some thoughts on this theorem from the perspective of recent results  on central limit theorem obtained via Malliavin calculus.

  \section{Preliminaries}

 Let $(W_{t})_{t\in [0,1]}$ be a classical
Wiener process on a standard Wiener space $\left( \Omega
,{\mathcal{F}},\mathbf{P}\right) $. If $f\in
L^{2}([0,1]^{n})$ with $n\geq 1$ integer, we introduce the multiple Wiener-It\^{o} integral of $f$ with respect to $W$. The basic references are the monographs \cite{Mal} or \cite{N}.
Let $f\in {\mathcal{S}_{n}}$ be an elementary functions with $n$
variables  that can be written as \begin{equation*}
f=\sum_{i_{1},\ldots ,i_{n}}c_{i_{1},\ldots ,i_{n}}1_{A_{i_{1}}\times
\ldots \times A_{i_{n}}}
\end{equation*}%
where the coefficients satisfy $c_{i_{1},\ldots ,i_{n}}=0$ if two indices $%
i_{k}$ and $i_{l}$ are equal and the sets $A_{i}\in
{\mathcal{B}}([0,1])$ are pairwise disjoint. For  such a step function $f$ we
define

\begin{equation*}
I_{n}(f)=\sum_{i_{1},\ldots ,i_{n}}c_{i_{1},\ldots
i_{n}}W(A_{i_{1}})\ldots W(A_{i_{n}})
\end{equation*}%
where we put $W(A)=\int_{0}^{1} 1_{A}(s)dW_{s}$. It can be seen that the application $%
I_{n}$ constructed above from ${\mathcal{S}}_{n}$ to $L^{2}(\Omega
)$ is an
isometry on ${\mathcal{S}}_{n}$ , i.e.%
\begin{equation}
\mathbf{E}\left[ I_{n}(f)I_{m}(g)\right] =n!\langle f,g\rangle
_{L^{2}([0,1]^{n})}\mbox{ if }m=n  \label{isom}
\end{equation}%
and%
\begin{equation*}
\mathbf{E}\left[ I_{n}(f)I_{m}(g)\right] =0\mbox{ if }m\not=n.
\end{equation*}%

Since the set ${\mathcal{S}_{n}}$ is dense in $L^{2}([0,1]^{n})$ for every $n\geq 1$ the mapping $%
I_{n}$ can be extended to an isometry from $L^{2}([0,1]^{n})$ to $%
L^{2}(\Omega)$ and the above properties hold true for this
extension.
It also holds that
\begin{equation}\label{ftilde}
I_{n}(f) = I_{n}\big( \tilde{f}\big)
\end{equation}
where $\tilde{f} $ denotes the symmetrization of $f$ defined by $\tilde{f}%
(x_{1}, \ldots , x_{n}) =\frac{1}{n!} \sum_{\sigma \in {\bf S} _{n}}
f(x_{\sigma (1) }, \ldots , x_{\sigma (n) } ) $.
We will need the general formula for calculating products of Wiener chaos
integrals of any orders $m,n$ for any symmetric integrands $f\in
L^{2}([0,1]^{\otimes m})$ and $g\in L^{2}([0,1]^{\otimes n})$; it is%
\begin{equation}
I_{m}(f)I_{n}(g)=\sum_{\ell=0}^{p\wedge q}\ell!C_{m}^{\ell}C_{n}^{\ell}%
I_{m+n-2\ell}(f\otimes_{\ell}g) \label{product}%
\end{equation}
where the contraction $f\otimes_{\ell}g$ is  defined by
\begin{eqnarray}
&&  (f\otimes_{\ell} g) ( s_{1}, \ldots, s_{m-\ell}, t_{1}, \ldots, t_{n-\ell
})\nonumber\\
&&  =\int_{[0,T] ^{m+n-2\ell} } f( s_{1}, \ldots, s_{m-\ell}, u_{1},
\ldots,u_{\ell})g(t_{1}, \ldots, t_{n-\ell},u_{1}, \ldots,u_{\ell})
du_{1}\ldots du_{\ell} . \label{contra}%
\end{eqnarray}
Note that the contraction $(f\otimes_{\ell} g) $ is an element of $L^{2}([0,1]^{m+n-2\ell})$ but it is not necessary symmetric. We will denote by $(f\tilde{\otimes }_{\ell} g)$ its symmetrization.

We recall that any square integrable random variable which is
measurable with respect to the $\sigma$-algebra generated by $W$ can
be expanded into an orthogonal sum of multiple stochastic integrals
\begin{equation}
\label{sum1} F=\sum_{n\geq0}I_{n}(f_{n})
\end{equation}
where $f_{n}\in L^{2}([0,1]^{n})$ are (uniquely determined)
symmetric functions and $I_{0}(f_{0})=\mathbf{E}\left[  F\right]  $.

We denote by $D$  the Malliavin  derivative operator that acts on smooth functionals of the form $F=g(W(\varphi _{1}), \ldots , W(\varphi_{n}))$ (here $g$ is a smooth function with compact support and $\varphi_{i} \in L^{2}([0,1])$ for $i=1,..,n$)
\begin{equation*}
DF=\sum_{i=1}^{n}\frac{\partial g}{\partial x_{i}}(W(\varphi _{1}), \ldots , W(\varphi_{n}))\varphi_{i}.
\end{equation*}
We can define the $i$ th Malliavin derivative $D^{(i)}$ is defined iteratively. The operator $D^{(i)}$ can be extended to the closure $\mathbb{D}^{p,2}$ of smooth functionals with respect to the norm
\begin{equation*}
\Vert F\Vert _{p,2}^{2} = \mathbf{E}F^{2}+ \sum_{i=1}^{p} \mathbf{E} \Vert D^{i} F\Vert ^{2} _{L^{2}([0,1] ^{i})}
\end{equation*}
 The adjoint of
$D$ is denoted by $\delta $ and is called the divergence (or
Skorohod) integral. Its domain $Dom(\delta)$  coincides with the class of stochastic processes $u\in L^{2}(\Omega \times [0,1])$ such that
\begin{equation*}
\left| \mathbf{E}\langle DF, u\rangle \right| \leq c\Vert F\Vert _{2}
\end{equation*}
for all $F\in \mathbb{D}^{1,2}$ and $\delta (u)$ is the element of $L^{2}(\Omega)$ characterized by the duality relationship
\begin{equation*}
\mathbf{E}(F\delta (u))= \mathbf{E}\langle DF, u\rangle.
\end{equation*}
For adapted integrands, the divergence integral coincides with
the classical It\^o integral.

Let $L$ be the Ornstein-Uhlenbeck operator defined on $Dom (L)= \mathbb{D}^{2,2}$
\begin{equation*}
LF=-\sum_{n\geq 0} nI_{n}(f_{n})
\end{equation*}
if $F$ is given by (\ref{sum1}). There exists a connection between $\delta, D $ and $L$ in the sense that a random variable $F$ belongs to the domain of $L$ if and only if $F\in \mathbb{D}^{1,2}$ and $DF \in Dom (\delta)$ and then $\delta DF=-LF$. Also we will need in the paper  the integration by parts formula
\begin{equation}
\label{ip}
F\delta(u)= \delta (Fu) + \langle DF, u\rangle
\end{equation}
whenever $F\in \mathbb{D}^{1,2}$, $u\in Dom(\delta)$ and $\mathbf{E}F^{2}\int_{0}^{1}u_{s}^{2}ds <\infty$.

\section{Asymptotic Cram\'er's theorem  }

We start by proving an asymptotic version of the Cram\'er's theorem \cite{Cramer}. A particular case (when the sequences $X_{n}$ and $Y_{n}$ are multiple integrals in a Wiener chaos of fixed order) has been proven in \cite{PT}, Corollary 1. Our proof is based on the Cram\'er's theorem (see \cite{Cramer}) and an idea from \cite{NOT}.

\begin{theorem}\label{cra}
Suppose that $(X_{n})_{n\geq 1}$ and $(Y_{n}) _{n\geq 1} $are two sequences of centered random variables in $L^{2}(\Omega)$ such that $\mathbf{E}X_{n} ^{2}\to _{n\to \infty }c_{1}>0$ and $\mathbf{E}Y_{n}^{2}\to _{n\to \infty} c_{2}>0$ with $c_{1}+ c_{2}=1$. Assume that for every $n\geq 1$, the random variables $X_{n}$ and $Y_{n}$ are independent. Then
\begin{equation*}
X_{n}+ Y_{n} \to N(0,1) \Leftrightarrow (X_{n}\to N(0, c_{1}) \mbox{ and } Y_{n} \to N(0, c_{2})).
\end{equation*}
 \end{theorem}
{\bf Proof: } One direction is  trivial. Let us assume that $X_{n} + Y_{n} \to _{n\to \infty} N(0,1)$. We will prove that $X_{n}\to _{n\to \infty} N(0, c_{1})$ and $Y_{n}\to _{n\to \infty} N(0, c_{2})$. Since $\mathbf{E}X_{n}^{2} \to_{n \to \infty}c_{1}$ and  $\mathbf{E}Y_{n}^{2}\to _{n\to \infty} c_{2}$ it follows that the sequence $(X_{n}, Y_{n})_{n\geq 1}$ is bounded in $L^{2}(\Omega)$. By Prohorov's theorem it suffices to prove that for any subsequence which converges  in distribution to some random vector $(F,G)$,  then we must have $F\sim N(0, c_{1}), G\sim N(0, c_{2}) $ and $F,G$ are independent.  Let us consider  such an arbitrary sequence $(X_{n_{k}}, Y_{n_{k}})$ which  converges in law  to $(F,G)$ as $k\to \infty$.  Because $X_{n_{k}}$ and $Y_{n_{k}}$ are independent for each $k$, it is clear that $F$ and $G$ are independent. Since $X_{n}+ Y_{n}\to _{n\to \infty} N(0,1)$ it follows that $F+G\sim N(0,1)$.

Cram\'er's theorem implies that $F\sim N(0, c_{1}) $ and $G\sim N(0, c_{2})$. \qed

\vskip0.3cm

This result can be extended to finite and even infinite sums of independent random variables.
\begin{prop}
Suppose that for every $n\geq 1$, $X^{n}= \sum_{k\geq 1} X^{n}_{k}$ where for every $n$ the random variables $X^{k}_{n}, k\geq 1$ are mutually independent  and the series is convergent for every $\omega$. Assume also that $X^{n}_{k}$ are centered for every $n,k \geq 1$ and $ \mathbf{E}(X^{n}_{k}) ^{2} \to _{n\to \infty} c_{k}>0$ for every $k\geq 1$. Suppose that $X^{n}$ converging in law to $N(0,1)$ as $n\to \infty$. Then for every $k\geq 1$ the sequence $X^{n}_{k}$ converges to the normal law as $n\to \infty$.
 \end{prop}
 {\bf Proof: } Since $X^{n} =X_{1}^{n}+ \sum_{k\geq 2} X_{k}^{n}$ and the two summands are independent, Theorem \ref{cra} implies that $X_{1}^{n}$ converges to the normal law. Inductively, the conclusion can be obtained. \qed

\begin{remark}
When, for every $n\geq 1$,  $X_{n}=I_{k_{1}}(f^{n}) $ and $Y_{n}=I_{k_{2}}(g^{n})$ are multiple stochastic integrals (possibly of different orders, that can also vary with $n$) we can go further by proving the following result. If $\mathbf{E} (X_{n} + Y_{n} ) ^{2} \to _{n\to \infty}1$ and $X_{n}+ Y_{n} $ converges in law to $N(0,1)$, if $\liminf _{n}\mathbf{E} X_{n}^{2} >0$ and $\liminf _{n}\mathbf{E} Y_{n}^{2} >0$
then
\begin{equation}
\label{rev1}
d_{Kol} (X_{n}, N(0, \mathbf{E}X_{n}^{2} )) \to _{n\to \infty }0 \mbox{ and } d_{Kol} Y_{n}, N(0, \mathbf{E}Y_{n}^{2} ) \to _{n\to \infty }0
\end{equation}
Here $d_{Kol} $ means the Kolmogorov distance (recall that the Kolmogorov distance between the law of two random variables $U$ and $V$ is given by $d_{Kol}(U,V)= \sup_{x\in \mathbb{R}}\left| P(U\leq x) -P(V\leq x)\right| $). That is, there is an asymptotic Cram\'er's theorem even if the variances of $X_{n}$ and $Y_{n}$ do not converge a priori. Relation (\ref{rev1}) can be proved as follows. First, recall the following bound when $X$ lives in a chaos of fixed order (see e.g. \cite{NoPe3})
\begin{equation}
\label{rev2}
d_{Kol} X, N(0, \mathbf{E}X^{2} ) \leq \frac{ \left( \left| \mathbf{E} X^{4} -3(\mathbf{E} X^{2})^{2} \right|\right)^{\frac{1}{2}} }{\mathbf{E}X^{2} } =: \frac{\left(| k_{4} (X) | \right) ^{\frac{1}{2}}}{\mathbf{E}X^{2} }
\end{equation}
where $k_{4}(X)$ is the fourth cumulant of $X$. It is immediate, by the definition of the cumulant, that $k_{4}(X+Y)=k_{4}(X)+ k_{4}(Y)$ if $X$ and $Y$ are independent. Moreover, it follows from \cite{NoPe3}, identity (3. 31) that $k_{4}(X)\geq 0$ if $X$ is a multiple integral. Hence, if $\mathbf{E} (X_{n} + Y_{n} ) ^{2} \to _{n\to \infty}1$ and $X_{n}+ Y_{n} $ converges in law to $N(0,1)$, then
$$k_{4}(X_{n})+ k_{4}(Y_{n})=k_{4}(X_{n}+ Y_{n}) =\mathbf{E}(X_{n}+Y_{n}) ^{4} -3 (\mathbf{E}(X_{n}+Y_{n})^{2})^{2} \to _{n\to \infty}0$$
and this implies that $k_{4}(X_{n})\to _{n\to \infty}0$ and $k_{4}(Y_{n}) \to_{n\to \infty}0$. The convergence (\ref{rev1}) is obtained by using (\ref{rev2}) and the hypothesis $\liminf _{n}\mathbf{E} X_{n}^{2} >0$ and $\liminf _{n}\mathbf{E} Y_{n}^{2} >0$.
\end{remark}

 \vskip0.3cm

\section{Decoupling central limit under strong independence}

Let us regard Theorem \ref{cra} from the perspective of the results in \cite{NoPe1}. In this part all random variables are centered. We  recall some facts related to the convergence of a sequence of random variables to the normal law in terms of the Malliavin calculus. For any random variable $X \in \mathbb{D}^{1,2}$ we  denote by
\begin{equation}\label{gx}
G_{X} := \langle DX, D(-L) ^{-1}X\rangle.
\end{equation}
The following result is a slight extension of Proposition 3.1 in \cite{NoPe1}. See also Theorem 3 in \cite{T1}.

\begin{prop}\label{t1}
Let $(X_{n})_{n\geq 1}$ be a sequence of square integrable random variables such that $\mathbf{E}X_{n}^{2} \to _{n\to \infty} c>0$. Then the following are equivalent:
\begin{description}
\item{1. }The sequence $(X_{n})_{n\geq 1}$ converges in law, an $n\to \infty$,  to the normal random variable $N(0,c)$, $c>0$.
\item{2. }For every $t\in \mathbb{R}$,
$\mathbf{E}\left( e^{itX_{n} } (c-G_{X_{n}})\right)\to _{n\to \infty }0.$
\item{3. }$
\mathbf{E}\left( (c-G_{X_{n}})
| X_{n}\right) \to _{n\to \infty} 0$ \mbox{ a.s. }.
\item{4. }For every $z\in \mathbb{R}$,
$\mathbf{E}\left(  f'_{z} (X_{n}) (c-G_{X_{n}})\right) \to_{n\to \infty }0.$
\end{description}
\end{prop}
{\bf Proof: } We follow the scheme $1. \Rightarrow  2. \Rightarrow 3.\Rightarrow 4. \Rightarrow 1$. The implications $1. \Rightarrow  2.$ and $3.\Rightarrow 4. \Rightarrow 1$ follow exactly as in \cite{T1}, Theorem 3. Concerning $2. \Rightarrow 3.$, set $F_{n}=c-G_{X_{n}}$ for every $n\geq 1$. The random variable $\mathbf{E}(F_{n}|X_{n})$ is the Radon-Nykodim derivative with respect to $P$ of the measure $Q_{n}(A)= \mathbf{E}(F_{n}1_{A})$, $A\in \sigma(X_{n})$. Relation 2. means that
$\mathbf{E}\left( e^{itX_{n}} \mathbf{E}(F_{n} / X_{n})\right)=\mathbf{E}_{Q_{n} } (e^{itX_{n}}) \to _{n\to \infty} 0$  and hence $\int_{\mathbb{R}}e^{ity}d(Q_{n}\circ X_{n}^{-1})(y)\to _{n\to \infty}0$. This implies that $Q_{n} (A)=\mathbf{E}(F_{n}1_{A}) \to _{n\to \infty} 0$  for any $A \in \sigma (X_{n})$ or $\mathbf{E}(F_{n}|X_{n}) \to _{n\to \infty}0.$
\vskip0.2cm
As  an immediate consequence we have (see also \cite{NoPe1}).

\begin{corollary}\label{np}
Suppose that $(X_{n})_{n\geq 1}$ is a sequence of random variables such that $\mathbf{E}X_{n}^{2} \to_{n \to \infty} c$. suppose that
\begin{equation}\label{s1}
\mathbf{E}(c-G_{X_{n}})^{2}\to _{n\to \infty }0.
\end{equation}
Then $X_{n} \to _{n\to \infty } N(0,c)$.
\end{corollary}

\begin{remark}\label{chafin}
In the case when the variables $X_{n}$ live in a fixed Wiener chaos, $X_{n}= I_{k}(f_{n})$, then the convergence in distribution of $X_{n}$ to the normal law is equivalent to (\ref{s1}), see \cite{NOT}.
\end{remark}
\vskip0.3cm

Assume that $(X_{n})_{n\ge 1}$ and $(Y_{n})_{n\geq 1}$ are two sequences of random variables such that: i) for every $n\geq 1$ the random variables $X_{n}$ and $Y_{n}$ are independent and ii) $X_{n}+ Y_{n}\to N(0,1)$ in distribution as $n\to \infty$.
 The quantity $G_{X_{n} + Y_{n}}$, which plays a central role, can be written as
\begin{equation*}
G_{X_{n} + Y_{n}}= G_{X_{n}}+ G_{Y_{n}}+ \langle DX_{n} , D(-L)^{-1}Y_{n}\rangle + \langle DY_{n} , D(-L)^{-1}X_{n}\rangle.
\end{equation*}
The force of the Cram\'er's theorem can be observed here: the fact that
\begin{equation*}
 \mathbf{E} \left( c_{1}-G_{X_{n}}+ c_{2}-G_{X_{n}} -\langle DX_{n} , D(-L)^{-1}Y_{n}\rangle - \langle DY_{n} , D(-L)^{-1}X_{n}\rangle| X_{n}+Y_{n}\right)
\end{equation*}
converges to zero implies that $\mathbf{E}(c_{1}-G_{X_{n}}| X_{n}) $ and $\mathbf{E}(c_{2}-G_{Y_{n}}| Y_{n}) $ both converge to zero. It is  not obvious to prove this  directly. Note also that the independence of $X_{n}$ and $Y_{n}$ does not guarantee a priori that the terms $\mathbf{E}( \langle DY_{n} , D(-L)^{-1}X_{n}\rangle |X_{n}+Y_{n})$ $\mathbf{E}( \langle DX_{n} , D(-L)^{-1}Y_{n}\rangle |X_{n}+Y_{n})$ vanish.
But the situation when these two terms vanish is also interesting and we will analyze this case in the sequel. We will see that it requires a slightly stronger assumption than the independence of $X_{n}$ and $Y_{n}$. We introduce the following concept.

\begin{definition} Two random variables $X=\sum_{n\geq 0} I_{n}(f_{n})$ and $Y=\sum_{m\geq 0} I_{m}(g_{m}) $ are called {\it strongly independent} if  for every $m,n\geq 0$, the random variables $I_{n}(f_{n})$ and $I_{m}(g_{m})$ are  independent.
\end{definition}

\begin{remark}\label{indep}
Let us recall the criterion for the independence of two multiple integrals given in \cite{UZ}:  Let $X'=I_{p}(f)$ and $Y'=I_{q}(g)$ where $f\in L^{2}([0,1] ^{p})$ and $g\in L^{2}([0,1] ^{q})$ ($p,q\geq 1$) are symmetric functions. Then $X'$ and $Y'$ are independent if and only if
\begin{equation*}
f\otimes _{1}g=0 \mbox{ almost everywhere on  }[0,1] ^{p+q-2}.
\end{equation*}
As a consequence two random variables $X$ and $Y$ as in Definition 1 are strongly independent if and only if for every $m,n\geq 1$,  $f_{n} \otimes _{1} g_{m}=0$ almost everywhere on $[0,1]^{m+n-2}$.

\end{remark}

Let us also note that the class of strongly independent random variables is strictly included in the class of independent random variables. Indeed, consider
\begin{equation*}
X_{1}=\sqrt{2} I_{1} (1_{[\frac{1}{2}, 1]}) \mbox{ and }Y_{1}= \sqrt{2} \int_{0}^{\frac{1}{2}} sign (W_{s}) dW_{s}.
\end{equation*}
Then $X_{1}$ and $Y_{1}$ are independent standard normal random variables. Define
\begin{equation*}
X=\frac{1}{\sqrt{2}}(X_{1}+ Y_{1} ) \mbox{ and } Y= \frac{1}{\sqrt{2}}(X_{1}- Y_{1} ).
\end{equation*}
Then $X, Y$ are also independent standard normal but they are not strongly independent because for example the chaoses of order one of $X$ and $Y$ are not independent (note that the random variable $ \int_{0}^{\frac{1}{2}} sign (W_{s}) dW_{s}$ has only even order chaos components).

\begin{lemma}\label{n11}
Assume that $X,Y\in \mathbb{D}^{1,2}$ and $X,Y$ are strongly independent. Then
\begin{equation*}
\langle DX, D(-L)^{-1}Y\rangle =0  \hskip0.5cm \mbox{ a.s }.
\end{equation*}
\end{lemma}
{\bf Proof: } Suppose first that $X=I_{n}(f)$ and $Y=I_{m}(g)$. Then, since $D_{\alpha } X=n I_{n-1}(f(\cdot , \alpha ))$ and $D_{\alpha}(-L)^{-1}Y= I_{m-1}(g(\cdot , \alpha))$,  using (\ref{product})
\begin{eqnarray*}
\langle DX, D(-L)^{-1}Y\rangle &=& n\int_{0}^{1} d\alpha I_{n-1}(f(\cdot, \alpha )) I_{m-1}(g(\cdot , \alpha ))\\
&=&m \sum_{k=0}^{(m-1)\wedge (n-1)} k! C_{m-1}^{k}C_{n-1}^{k} \int_{0}^{1}d\alpha I_{m+n-2-2k} (f(\cdot , \alpha )\otimes _{k} g(\cdot , \alpha ))\\&=&
\sum_{k=0}^{(m-1)\wedge (n-1)} k! C_{m-1}^{k}C_{n-1}^{k} I_{m+n-2-2k}(f\otimes _{k+1} g)
\end{eqnarray*}
and this is equal to zero from the characterization of the independence of two multiple multiple integrals (see Remark \ref{indep}). The extension to the general case is immediate since, if $X=\sum_{n}I_{n}(f_{n})$ and $Y=\sum_{m}I_{m}(g_{m})$,
\begin{equation*}
\langle DX, D(-L) ^{-1}Y\rangle =\sum_{m,n}\langle DI_{n}(f_{n}), D(-L)^{-1}I_{m}(g_{m})\rangle .
\end{equation*} \qed

\vskip0.3cm

In view of Lemma \ref{n11}, the Proposition \ref{t1} can be formulated for strongly independent random variables as follows: Suppose that $(X_{n})_{n\geq 1}$ and $(Y_{n})_{n\geq 1}$ are two sequences of centered strongly independent random variables such that $\mathbf{E}X_{n}^{2}\to _{n\to \infty }c_{1} $ and $\mathbf{E}Y_{n}^{2} \to _{n \to \infty }c_{2}$ where $c_{1}, c_{2}>0$ are such that $c_{1}+ c_{2} =1$. Then the following affirmations are equivalent:
 \begin{description}
 \item{1. }The sequence $(X_{n}+ Y_{n})_{n\geq 1}$ converges in law to a standard normal random variable as $n\to \infty$;
 \item{2. } For every $t\in \mathbb{R}$, $\mathbf{E}\left( e^{it(X_{n}+Y_{n}) }(c_{1}-G_{X_{n}}+ c_{2}-G_{X_{n}}) \right) \to _{n\to \infty }0.$;
 \item{3. } $\mathbf{E}\left( c_{1}-G_{X_{n}}+ c_{2}-G_{Y_{n}} | X_{n}+ Y_{n}\right) \to _{n\to \infty}0.$;
  \item{4. } For every $z\in \mathbb{R}$, $\mathbf{E}\left( f'_{z}(X_{n}+ Y_{n})(c_{1}-G_{X_{n}}+ c_{2}-G_{Y_{n}}) \right) \to_{n\to \infty }0$.
\end{description}

\vskip0.3cm

Let us assume now that the two sequences of Theorem \ref{cra} are strongly independent. We will also  assume that the convergence of the sum $X_{n}+ Y_{n}$ to $N(0,1)$ is strong in the sense that $\mathbf{E}\left( 1-G_{X_{n}+Y_{n}}\right) ^{2}$ converges to zero as $n\to \infty$. We can say, somehow,  that in this case the sum $X_{n}+Y_{n}$ is rather close to the normal law since the upper bound of the distance between it and $N(0,1)$ goes to zero. We will prove that this implies that the convergence of $X_{n}$ and $Y_{n}$ to the normal law is also strong.

\begin{remark}
The case of multiple stochastic integrals can be easily understood. Suppose that $X_{n}=I_{k}(f^{n})$ and $Y_{n}=I_{l}(g^{n}) $ where for every $n\geq 1$ the kernels $f^{n}, g^{n}$ are in $L^{2}([0,1] ^{k})$ and $L^{2}([0,1]^{l})$ respectively. Assume that $EX_{n}^{2}\to_{n \to \infty}c_{1}>0$ and $EY_{n}^{2}\to_{n \to \infty}c_{2}>0$ such that $c_{1}+c_{2}=1$. Then if $X_{n}+ Y_{n}\to_{n \to \infty} N(0,1)$ and $X_{n}, Y_{n}$ are independent (thus strongly independent) it follows that $X_{n}\to N(0,c_{1})$ and $Y_{n}\to N(0,c_{2})$ and by Remark \ref{chafin}, $E\left( c_{1}-G_{X_{n}}\right) ^{2} \to_{n \to \infty}0$ and $E\left( c_{2}-G_{Y_{n}}\right) ^{2} \to_{n \to \infty}0$, so the convergence of $X_{n}$ and $Y_{n}$  to the normal distribution is strong.
\end{remark}
\vskip0.2cm

We will also need the following lemma.

\begin{lemma}\label{n12}
Assume that $X,Y\in \mathbb{D}^{1,2}$ and $X,Y$ are strongly independent. Then  the random variables $G_{X}$ and $G_{Y}$ are strongly independent.
\end{lemma}
{\bf Proof: } Let us assume once again that $X=I_{n}(f)$ and $Y=I_{m}(g)$. The result can  easily be extended to the general case.
We have
\begin{equation*}
G_{X}=n \sum_{k=0}^{n-1} \left( C_{n-1}^{k} \right) ^{2} I_{2n-2-2k}(f\otimes _{k+1} f)
\end{equation*}
and
\begin{equation*}
G_{Y}= m\sum_{l=0}^{m-1} l! \left( C_{m-1}^{l}\right) ^{2} I_{2m -2-2l} (g\otimes _{l+1} g).
\end{equation*}
It suffices to show that for every $k=1,..,n-1$ and $l=1, .., m-1$ the random variables $I_{2n-2k}(f\otimes _{k} f) $ and $I_{2m-2l} (g\otimes _{l}g)$ are independent or equivalently
\begin{equation*}
(f\tilde{\otimes }_{k}f) \otimes _{1} (g\tilde{\otimes }_{l} g)= 0 \mbox{ a.e.  on } [0,1] ^{2m-2k+2m-2l-2}.
\end{equation*}
But since
\begin{eqnarray*}
&&(f\tilde{\otimes }_{k}f)(x_{1},..,x_{2n-2k})\\
 &&=\sum_{\sigma \in {\bf S}_{2n-2k} } \int_{[0,1] ^{k}} f(u_{1},..,u_{k}, x_{\sigma (1)},..,x_{\sigma (n-k)})f(u_{1},..,u_{k}, x_{\sigma (n-k+1)},..,x_{\sigma (2n-2k)})du_{1}..du_{k}
\end{eqnarray*}
and
\begin{eqnarray*}
&&(g\tilde{\otimes }_{l} g)(y_{1},..,y_{2m-2l})\\
&& \sum_{\rho \in {\bf S}_{2m-2l}} \int_{[0,1]^{l}}g(v_{1},..,v_{l}, y_{\rho(1)},..,y_{\rho (m-l) })g(v_{1},..,v_{l}, y_{\rho(m-l+1)},..,y_{\rho (2m-2l) })
dv_{1}..dv_{l}
\end{eqnarray*}
then $(f\tilde{\otimes }_{k}f) \otimes _{1} (g\tilde{\otimes }_{l} g)= 0 $ almost everywhere on $[0,1] ^{2m-2k+2m-2l-2}$ by using Fubini and the fact that $\int_{0}^{1} d\alpha f(u_{1},..,u_{k}, x_{1},..,x_{n-k-1}, \alpha ) g(v_{1},..,v_{l}, y_{1}, ..,y_{m-l+1}, \alpha  ) =0$ for almost all $u_{1},v_{i}, x_{i},y_{i}$.
The general case demands to prove that $(f_{n} \tilde{\otimes} _{k}f_{n'} ) \otimes _{1} (g_{m}\tilde{\otimes }_{l} g_{m'})=0$ almost everywhere for every $n,n',m,m'\geq 1$ and for every  $k=1,..,n\wedge n '$ and $l=1, .., m\wedge m'$ and this can be done similarly as above (note that the fact that $k,l\geq 1$ and the value zero is excluded is essential for the proof). \qed

\begin{prop}\label{p101}
Suppose that $(X_{n})_{n\geq 1}$ and $(Y_{n})_{n\geq 1}$ are two sequences of centered strongly independent random variables such that $\mathbf{E}X_{n}^{2}\to _{n\to \infty }c_{1} $ and $\mathbf{E}Y_{n}^{2} \to _{n \to \infty }c_{2}$ where $c_{1}, c_{2}>0$ are such that $c_{1}+ c_{2} =1$. Then $\mathbf{E}\left(1-G_{X_{n}+Y_{n}}\right) ^{2} \to _{n\to \infty }0 $ if and only if
\begin{equation*}
\mathbf{E}\left( c_{1}-G_{X_{n}}\right) ^{2} \to _{n\to \infty }0 \mbox{ and }  \mathbf{E}\left( c_{2}-G_{Y_{n}}\right) ^{2} \to _{n\to \infty }0.
\end{equation*}
\end{prop}
{\bf Proof: } By using Lemmas \ref{n11} and \ref{n12} we have
\begin{equation*}
\mathbf{E}\left(1-G_{X_{n}+Y_{n}}\right) ^{2}=\mathbf{E}\left( c_{1}-G_{X_{n}}\right) ^{2}+\mathbf{E}\left( c_{2}-G_{Y_{n}}\right) ^{2}
\end{equation*}
and the conclusion is immediate. \qed

\vskip0.3cm

We will study  next if the result in Proposition \ref{p101} can be obtained by relaxing the hypothesis on the strong independence of $X_{n}$ and $Y_{n}$ for every $n$. As we have seen, the strong independence of two variables $X$ and $Y$ implies that $\langle DX, D(-L)^{-1}Y\rangle = \langle DY, D(-L)^{-1}X\rangle =0$ a.s. But in order to eliminate the "mixed" terms we only need $\mathbf{E}( \langle DX, D(-L)^{-1}Y\rangle |X+Y) =\mathbf{E}( \langle DY, D(-L)^{-1}X\rangle | X+Y) =0$ a.s. We  therefore introduce an intermediary class between the class of independent random variables and the class of strongly independent random variables for which this property holds.

\begin{definition}
We will say that a couple $(X,Y)$  of two random variables  in the space $\mathbb{D}^{1,2}$ belongs to the class ${\cal{A}}$  if the vector $X$ is independent of the vector $(Y, G_{Y})$ and $Y$ is independent of the vector  $(X, G_{X})$.
\end{definition}

We will give now examples of random variables in ${\cal{A}}$. First we recall the following result from \cite{UZ}.
\begin{lemma}\label{uz}
Let $X\in \mathbb{D}^{1,2}$ and $Y,Z\in L^{2}(\Omega)$. Then $X $ is independent of the pair $(Y,Z)$ if and only if for every $\alpha , \beta \in \mathbb{R}$
\begin{equation*}
\mathbf{E}\left( \langle DX, D(-L)^{-1}e^{i(\alpha Z + \beta Y)}\rangle | X \right)= 0 \mbox{ a.s. }.
\end{equation*}
\end{lemma}

We show that a couple of strongly independent random variables is in the set ${\cal{A}}$. We consider first the case of multiple integrals.
\begin{lemma}
Suppose that $X=I_{p}(f)$ and $Y=I_{q}(g)$ where $f\in L^{2}([0,1] ^{p}) $ and $g\in L^{2}([0,1] ^{q} )$ ($p,q\geq 1$)  are symmetric functions. Assume that $X$ and $Y$ are independent. Then $(X, Y)$ belongs to the class ${\cal{A}}$.
\end{lemma}
{\bf Proof: }We will prove that $X$ is independent of the couple $(Y, G_{Y})$. Similarly it will follow that $Y$ is independent of $(X, G_{X})$. We prove that
\begin{equation*}
\langle DX, D(-L)^{-1}e^{i(\alpha Y + \beta G_{Y}) } \rangle =0 \mbox{ a.s. }
\end{equation*}
or, since $D(-L) ^{-1}= (I+L)^{-1}D $,
\begin{equation*}
\langle DX, (I+L)^{-1}De^{i(\alpha Y + \beta G_{Y}) } \rangle =0  \mbox{ a.s. }
\end{equation*}
Note that $De^{i(\alpha Y + \beta G_{Y}) }= e^{i(\alpha Y+ \beta G_{Y})} (i\alpha DY + i\beta DG_{Y})$. First we will show that
\begin{equation*}
\langle DX, e^{i(\alpha Y+ \beta G_{Y})}DY \rangle =0 \mbox{ a. s. }
\end{equation*}
Assume that the random variable $e^{i(\alpha Y+ \beta G_{Y})}$ admits the chaos expansion $e^{i(\alpha Y+ \beta G_{Y})}=\sum_{N\geq 0} I_{N}(h_{N})$ (in the sense that its real part and its imaginary part admit  such a decomposition). Then
\begin{eqnarray*}
e^{i(\alpha Y+ \beta G_{Y})}D_{\alpha }Y &=& q\sum_{N\geq 0} I_{N}(h_{N}) qI_{q-1}(g(\cdot , \alpha ))\\
&=& q\sum_{N\geq 0} \sum_{r=0} ^{N\wedge (q-1)} r! C_{q-1}^{r} C_{N}^{r} I_{N+q-1-2r} (h_{N}\otimes _{r} g(\cdot , \alpha ))
\end{eqnarray*}
and
\begin{eqnarray*}
(I+L)^{-1}e^{i(\alpha Y+ \beta G_{Y})}D_{\alpha }Y &=&q\sum_{N\geq 0} \sum_{r=0} ^{N\wedge (q-1)} r! C_{q-1}^{r} C_{N}^{r} (1+N+(q-1)-2r)^{-1}I_{N+q-1-2r} (h_{N}\otimes _{r} g(\cdot , \alpha )).
\end{eqnarray*}
Therefore
\begin{eqnarray*}
&& \langle DX, (I+L)^{-1}De^{i(\alpha Y + \beta G_{Y}) } \rangle \\
&=&pq\sum_{N\geq 0} \sum_{r=0} ^{N\wedge (q-1)} r! C_{q-1}^{r} C_{N}^{r} (1+N+(q-1)-2r)^{-1}\\
&&\times \sum_{a=0}^{(N+q-1-2r)\wedge (p-1)}I_{N+q-1-2r+p-2a}\int_{0}^{1}
\left( (h_{N}\tilde{ \otimes } _{r} g(\cdot , \alpha )\otimes _{a} f(\cdot , \alpha )\right) d\alpha.
\end{eqnarray*}
Above, $(h_{N}\tilde{ \otimes } _{r} g(\cdot , \alpha )$ means the symmetrization of the function $(t_{1}, \ldots , t_{N+q-1-2r}) \to (h_{N} \otimes  _{r} g(t_{1}, \ldots , t_{N+q-1-2r} , \alpha )$ for fixed $\alpha$. In other words the above symmetrization does not affect the variable $\alpha$.  By interchanging the order of integration to integrate first with respect to $\alpha$ we will obtain that the last quantity is zero. Similarly it will follow that $\langle DX, e^{i(\alpha Y+ \beta G_{Y})}DG_{Y}\rangle $ is almost surely zero.
\qed

\vskip0.2cm

We can extend the previous result  to the case of infinite chaos expansion.
\begin{lemma}
Assume that $X,Y$ are two strongly independent random variables in $\mathbb{D}^{1,2}$.
Then $(X, Y)$ belongs to ${\cal{A}}$.
\end{lemma}
{\bf Proof: }The proof follows the lines of the proof of Lemma 4. In order to check that
\begin{equation*}
\langle DX, (I+L)^{-1}De^{i(\alpha Y + \beta G_{Y}) } \rangle =0 \mbox{ a.s. }
\end{equation*}
we write
\begin{eqnarray*}
&& \langle DX, (I+L)^{-1}De^{i(\alpha Y + \beta G_{Y}) } \rangle \\
&=&\sum_{p,q\geq 1}pq\sum_{N\geq 0} \sum_{r=0} ^{N\wedge (q-1)} r! C_{q-1}^{r} C_{N}^{r} (1+N+(q-1)-2r)^{-1}\\
&&\times \sum_{a=0}^{(N+q-1-2r)\wedge (p-1)}I_{N+q-1-2r+p-2a}\int_{0}^{1}
\left( (h_{N}\tilde{\otimes } _{r} g(\cdot , \alpha ))\otimes _{a} f(\cdot , \alpha )\right) d\alpha
\end{eqnarray*}
and we can finish as in the proof of the previous lemma. \qed

\vskip0.3cm

An interesting property of the couples in ${\cal{A}}$ is that the conditional expectation given $X+Y$ of the mixed scalar products $\langle DX, D(-L)^{-1} Y\rangle$ and $\langle DY, D(-L)^{-1} X\rangle$ vanish.

\begin{lemma}\label{a1}
Assume that $(X,Y)$  belongs to the class ${\cal{A}}$. Then
\begin{equation*}
\mathbf{E}\left(   e^{it (X+Y)}  \langle DX, D(-L)^{-1} Y\rangle \right) =0 \mbox{ a.s. }
\end{equation*}
\end{lemma}
{\bf Proof: } We have
\begin{eqnarray*}
\mathbf{E}\left(   e^{it (X+Y)}  \langle DX, D(-L)^{-1} Y\rangle \right) &=& \mathbf{E} \frac{1}{it}\langle  De^{itX}, e^{itY}D(-L)^{-1}Y \rangle \\
&=& \frac{1}{it} \mathbf{E} \left(  e^{itX}  \delta (e^{itY} D(-L)^{-1}Y) \right)\\
&=&\frac{1}{it} \mathbf{E}\left( e^{itX} \left( e^{itY}\delta D(-L)^{-1}Y   -it  e^{itY}\langle DY, D(-L)^{-1}Y\rangle \right)\right) \\
&=&\frac{1}{it} \mathbf{E}\left( e^{itX} \left( e^{itY}Y   -it  e^{itY}G_{Y} \right)\right)
\end{eqnarray*}
where we used the fact that since $e^{itY}\in \mathbb{D}^{1,2}$ and $D(-L)^{-1}Y \in Dom (\delta)$, then $e^{itY}D(-L^{-1})Y \in Dom (\delta )$ and by (\ref{ip})
$\delta (e^{itY}(D(-L)^{-1}Y))= e^{itY}\delta (D(-L)^{-1}Y) -ite^{itY}\langle DY, D(-L)^{-1}Y\rangle$.
By using the fact that $(X,Y)$ belongs to the class ${\cal{A}}$  we obtain
\begin{eqnarray*}
\mathbf{E}\left(   e^{it (X+Y)}  \langle DX, D(-L)^{-1} Y\rangle \right) &=&\frac{1}{it} \mathbf{E}(e^{itX}) \left(   \mathbf{E}(e^{itY}Y ) -it\mathbf{E}( e^{itY}G_{Y}) \right).
\end{eqnarray*}
 Now, going in the converse direction
 \begin{equation*}
 \mathbf{E}(e^{itY}Y ) -it\mathbf{E}( e^{itY}G_{Y}) =\mathbf{E}\delta (e^{itY}(D(-L)^{-1}Y))=0.
\end{equation*} \qed

We are now answering the following question: let $(X_{n})_{n\geq 1}$ and $(Y_{n})_{n\geq 1}$ be two sequences of random variables such that for every $n\geq 1$ the couple $(X_{n},Y_{n})$ is  in the class ${\cal{A}}$. Suppose that the sum $X_{n}+ Y_{n}$ converges  to the normal law and is such that the upper bound from Stein's method is attained, in the sense that $\mathbf{E}(1-G_{X_{n}+Y_{n} }) ^{2}$ converges to zero. Could we then conclude that both $X_{n}$ and $Y_{n}$ converge in a strong sense to the normal laws $N(0,c_{1}) $ and  $N(0, c_{2})$ respectively? We will see that this is true in some particular case under a supplementary hypothesis on the sequences $X_{n}$ and $Y_{n}$.

\subsection{Wiener chaos stable random variables}

Let us denote another class of families of random variables where the central limit of the sum implies central limit for each component. The idea is to assume a property on the filtration generated by $X_{n}+ Y_{n}$. Let us denote by $J_{n}$ the orthogonal projection of $L^{2}(\Omega)$ on the $n$-th Wiener chaos. We recall the following definition (see \cite{NUZ}, \cite{UZ2}).
\begin{definition}
We will say that a sigma-algebra $\tau \subset {\cal{F}}$ is Wiener chaos stable if for every $n$, $J_{n}\left( L^{2}(\tau)\right) \subset L^{2}(\tau)$. In other words, if a random variable $F\in L^{2}(\tau)$ admits the chaos decomposition  $F=\sum_{n\geq 0} I_{n}(f_{n})$ then for every $n\geq 0$ the random variable $I_{n}(f_{n})$ is $\tau$ -measurable.
\end{definition}

\begin{remark}
The Wiener chaos stable property for sigma-algebras is equivalent to the $L^{-1}$- stable property. Recall that a sigma-algebra $\tau $ is $L^{-1}$ stable if  $L^{-1}(L_{0}^{2}(\tau))\subset L^{2}_{0}(\tau)$ where $L^{2}_{0}(\tau)$ is the set of $\tau$-measurable square integrable random variables with zero expectation. As a matter of fact, the sigma -algebra generated by \\
$I_{p}(f), \langle DI_{p}(f), h_{1}\rangle, \langle DI_{p}(f), h_{2} \rangle, \ldots, \langle D^{p-1} I_{p}(f),h_{i_{1}}\otimes ..\otimes h_{i_{p-1}}\rangle $, where $h_{i}, i\geq 1$ is a complete orthogonal sequence in $L^{2}[0,1]$, is Wiener stable (see \cite{NUZ}, \cite{UZ2}).
\end{remark}

\begin{theorem}
Suppose that for every $X_{n}=\sum_{n\geq 1} I_{k}(f_{k}^{n})$ and $Y_{n}= \sum_{l\geq 1} I_{l}(g_{l}^{n}) $ are such that $\mathbf{E}X_{n}^{2} \to _{n\to \infty} c_{1}$ and $\mathbf{E}Y_{n}^{2}\to _{n\to \infty} c_{2}$ (such that $c_{1}, c_{2}>0$ and $c_{1}+ c_{2}=1$). Assume that
 \begin{description}
 \item{i. }for every $n\geq 1$ the couple $(X_{n}, Y_{n})$ belongs to the class ${\cal{A}}$.
 \item{ii. } for every $n\geq 1$ the sigma-algebras $\sigma (X_{n}) $ and $\sigma (Y_{n}) $ are Wiener chaos stable.
 \end{description}
Assume also that $\mathbf{E}\left( 1-G_{X_{n}+Y_{n}} \right) ^{2} \to _{n\to \infty }0.$ Then
\begin{equation*}
\mathbf{E}\left( c_{1}-G_{X_{n}}\right) ^{2} \to _{n\to \infty }0 \mbox{ and  } \mathbf{E}\left( c_{2}-G_{Y_{n}}\right) ^{2} \to _{n\to \infty }0.
\end{equation*}

\end{theorem}
{\bf Proof: } We will show that under assumption ii., the random variable $G_{X_{n}}$ belongs to $\sigma (X_{n})$ for every $n\geq 1$. Since $X_{n}$ is $\sigma(X_{n})$ measurable and $\sigma (X_{n}) $ is Wiener chaos stable, we get that $I_{k}(f_{k}^{n})$ is $\sigma(X_{n})$ measurable for every $n,k$. Consequently, $I_{k}^{n}(f_{k}^{n})I_{l}^{n}(f_{l}^{n})$ is $\sigma(X_{n})$ measurable for every $n,k,l$ and by using the product formula we will have that
\begin{equation*}
I_{k+l-2r} \left( f^{n}_{k} \otimes _{r} f^{n} _{l} \right)
\end{equation*}
is $\sigma(X_{n})$ measurable for every $n,k,l\geq 1$ and $r=0,.., k\wedge l $. As a consequence we can easily obtain that $G_{X_{n}}$ is measurable with respect to $\sigma(X_{n})$ and similarly $G_{Y_{n}}$ is measurable with respect to $\sigma(Y_{n})$.
Assume now that $\mathbf{E}(1-G_{X_{n}+Y_{n}})^{2} \to _{n\to \infty} 0$. The asymptotic Cram\'er's Theorem \ref{cra} together with Proposition \ref{t1} imply that $\mathbf{E}(c_{1}-G_{X_{n}}| X_{n}) \to 0$ and $\mathbf{E}(c_{2}-G_{Y_{n}}| Y_{n}) \to 0$  a.s. and by the measurability of $G_{X_{n}}$ and $G_{Y_{n}}$ we obtain the conclusion. \qed

\subsection{Vectorial convergence of $X_{n}+ Y_{n}$ and $G_{X_{n}}+ G_{Y_{n}}$}

A second class of sequences of random variables for which the central limit can be broken in order to ensure the strong convergence of each summand  is inspired by \cite{NoPe2}.
\begin{theorem}
Assume that $\mathbf{E}X_{n}^{2} \to c_{1}$ and $\mathbf{E}Y_{n}^{2}\to c_{2}$ (such that $c_{1}, c_{2}>0$ and $c_{1}+ c_{2}=1$). Assume that for every $n\geq 1$ the couple of random variables $(X_{n}, Y_{n})$ belongs to ${\cal{A}}$. Suppose moreover that the vector
\begin{equation}\label{cond3}
\left( X_{n}+ Y_{n}, \frac{c_{1}-G_{X_{n}}+ c_{2}-G_{Y_{n}}}{ \mathbf{E} \left( (c_{1}-G_{X_{n}})^{2} + (c_{2}-G_{Y_{n}}) ^{2} \right) ^{\frac{1}{2}}} \right)
\end{equation}
converges as $n\to \infty$ to the vector $(N_{1}, N_{2})$ where $N_{1}, N_{2}$ are standard normal random variables with correlation $\rho$. Then if $X_{n}+ Y_{n} \to_{n \to \infty} N(0,1)$ implies that
\begin{equation*}
\mathbf{E}\left(c_{1}-G_{X_{n}}\right) ^{2} \to_{n \to \infty}0 \mbox{ and } \mathbf{E}\left(c_{2}-G_{Y_{n}}\right) ^{2} \to_{n \to \infty}0
\end{equation*}
\end{theorem}
{\bf Proof: }On one hand, we have that
\begin{equation}\label{u1}
\mathbf{E}\left( f'_{z} (X_{n}+ Y_{n})(c_{1}-G_{X_{n}}+ c_{2}- G_{Y_{n}})\right) \to _{n\to \infty}0.
\end{equation}
On the other hand, from the convergence of the vector (\ref{cond3}) we get
\begin{equation*}
\mathbf{E}\left( f'_{z} (X_{n}+ Y_{n})(c_{1}-G_{X_{n}}+ c_{2}- G_{Y_{n}})a_{n}^{-1} \right) \to f'_{z}(N_{1})N_{2}
\end{equation*}
where $a_{n}= \mathbf{E} \left( (c_{1}-G_{X_{n}})^{2} + (c_{2}-G_{Y_{n}}) ^{2}  \right) ^{\frac{1}{2}}$.
It follows from the proof of Theorem 3.1 in \cite{NoPe2} that we can  find a constant $c\in (0,1)$ such that
\begin{equation}\label{u2}
\left| \mathbf{E}\left( f'_{z} (X_{n}+ Y_{n})(c_{1}-G_{X_{n}}+ c_{2}- G_{Y_{n}})\right)\right| ^{2} \geq c  \mathbf{E} \left( (c_{1}-G_{X_{n}})^{2} + (c_{2}-G_{Y_{n}}) ^{2} \right) ^{2}.
\end{equation}
By combining the  relations (\ref{u1}) and (\ref{u2}), we obtain that
$ \mathbf{E} \left( (c_{1}-G_{X_{n}})^{2} + (c_{2}-G_{Y_{n}}) ^{2} \right) ^{2}\to _{n\to \infty} 0 $ and this gives the convergence of $X_{n}$ and $Y_{n}$ to
$N(0, c_{1})$ and $N(0, c_{2})$ respectively. \qed

\subsection{Random variables with independent chaos components}

In this part we prove that in the case when the chaotic  components  appearing in the decomposition of $X_{n}$ are mutually independent (and  the same is true for  $Y_{n}$) then the central limit of the sum implies the central limit of the summands (in a strong sense) under simple independence.
\begin{prop}
Assume that for every $n\geq 1$, $X_{n} =\sum_{k\geq 1} I_{k}(f^{n}_{k})$ and $Y_{n}= \sum_{l\geq 1}I_{l} (g^{n}_{l}) $ and
 \begin{equation*}
 \mathbf{E}X_{n}^{2} \to_{n \to \infty} c_{1}, \hskip0.5cm \mathbf{E} Y_{n}^{2} \to_{n \to \infty} c_{2}
 \end{equation*}
 with  $c_{1}, c_{2}>0$ and $c_{1}+ c_{2}=1$.
 Suppose that the following conditions are fulfilled
 \begin{description}
\item{i. } for every $n\geq 1$ the random variables $X_{n}$ and $Y_{n}$ are  independent
\item{ii. } for every $n\geq 1$, the random variables $(I_{k}(f^{n}_{k}) )_{k\geq 1}$ are pairwise independent the same holds for $(I_{k} (g^{n}_{l}))_{l\geq 1}$.
\end{description}
Then $X_{n}+Y_{n} \to N(0,1)$ implies
\begin{equation*}
\mathbf{E}\left( c_{1}-G_{X_{n}}\right) ^{2} \to_{n \to \infty}0 \mbox{ and } \mathbf{E}\left( c_{2}-G_{Y_{n}}\right) ^{2} \to_{n \to \infty}0.
\end{equation*}
\end{prop}
{\bf Proof: } The Theorem \ref{cra} implies that $X_{n}\to N(0, c_{1})$ and $Y_{n}\to N(0,1)$ in law. Corollary  1 and Assumption ii. gives that for every $k$ the sequence $I_{k}(f_{l}^{n})$ converges to a normal law as $k\to \infty$. Finally, we use  Remark  \ref{chafin} and Lemmas \ref{n11}, \ref{n12}. \qed

{\bf Acknowledgments. }We wish to thank the anonymous referee for valuable  comments on our manuscript.


\begin{thebibliography}{99}

\bibitem{Cramer}
{H. Cram\'er (1936): } \"Uber eine Eigenschaft der normalen Verteilungsfunction. \emph{Math. Z., } 41(2), 405-414.

\bibitem{Mal}
{P.  Malliavin, P. (2002): }\emph{ Stochastic Analysis.} Springer-Verlag.

\bibitem {NoPe1}{I. Nourdin and G. Peccati (2009): }Stein's method on Wiener
chaos.  \emph{Probability Theory and Related Fields, }145 (1), 75-118.

\bibitem{NoPe2}
{I. Nourdin and G. Peccati (2009): } Stein's method and exact Berry -Ess\'{e}en asymptotics for functionals of Gaussian fields.  \emph{The Annals of Probability, }37(6), 2200-2230.

\bibitem{NoPe3}
{I. Nourdin and G. Peccati (2008): } Stein's method meets Malliavin calculus: a short survey with new estimates. To appear in \emph{Recent Advances in Stochastic Dynamics and Stochastic Analysis, } World Scientific.

\bibitem{N} {D. Nualart (2006): }\emph{Malliavin Calculus and Related
Topics. Second Edition. }{Springer. }

\bibitem {NOT}{D. Nualart and S. Ortiz-Latorre (2008): }Central limit theorems
for multiple stochastic integrals and Malliavin calculus. {\emph{Stochastic
Processes and their Applications, }\textbf{118}, 614-628.}

\bibitem{NUZ}
{D. Nualart, A.S. \"Ustunel and M. Zakai (1990): }Some relations among classes of $\sigma$ fields on Wiener space. \emph{Probability Theory and Related Fields, }84, 119-129.

\bibitem {PT}{G. Peccati and C.A. Tudor (2004): }Gaussian limits for
vector-valued multiple stochastic integrals. \emph{S\'{e}minaire de
Probabilit\'{e}s}{, \textbf{XXXIV}, 247-262.}

\bibitem{T1}
{C.A. Tudor (2009): }{\em On the structure of Gaussian random variables. } Preprint.

\bibitem{UZ}
{A.S. Ustunel and M. Zakai (1989): }On independence and conditioning on Wiener space. \emph{The Annals of Probability, } 17(4), 1441-1453.

\bibitem{UZ2} {A.S. Ustunel and M. Zakai (1990): }On the structure of independence on Wiener space. \emph{Journal of Functional Analysis, } 90, 113-137.

\end{thebibliography}
\end{document}